\documentclass{amsart}

\usepackage{fancyhdr}
\usepackage{lastpage}
\usepackage{stmaryrd,yhmath}

\newcommand{\bdis}{\begin{displaymath}}
\newcommand{\edis}{\end{displaymath}}
\newcommand{\be}{\begin{equation}}
\newcommand{\ee}{\end{equation}}

\newcommand{\mbb}{\mathbb}
\newcommand{\mcal}{\mathcal}

\newcommand{\vp}{\varphi}

\newcommand{\vth}{\vartheta}

\newcommand{\zf}{\zeta\left(\frac{1}{2}+it\right)}

\newtheorem{lemma}[]{Lemma}

\theoremstyle{definition}

\newtheorem{cor}[]{Corollary}

\theoremstyle{remark}
\newtheorem{remark}[]{Remark}

\newtheorem*{mydef1}{{\bf Theorem}}

\newtheorem*{mydef61}{{\bf Formula 1}}

\newtheorem*{mydef62}{{\bf Formula 2}}

\numberwithin{equation}{section}

\begin{document}

\title{The validity of the analog of the Riemann hypothesis for some parts of $\zeta(s)$ and the new formula for $\pi(x)$}

\author{Jan Moser}

\address{Department of Mathematical Analysis and Numerical Mathematics, Comenius University, Mlynska Dolina M105, 842 48 Bratislava, SLOVAKIA}

\email{jan.mozer@fmph.uniba.sk}


\begin{abstract}
An analog of the Riemann hypothesis is proved in this paper. Some new integral equations for the functions $\pi(x)$ and $R(x)$ follows. A new
effect that is shown is that these function - with essentially different behavior - are the solutions of the similar integral equations. \\

\noindent This paper is the English version of the paper of reference \cite{1}.
\end{abstract}

\maketitle

\section{The main result}

\subsection{}

Let (comp. \cite{2}, (7), (21); $2P\beta<\ln P_0$)

\be \label{1.1}
P=(\ln P_0)^{1-\epsilon},\ \beta=\left[\ln^{\frac{2\epsilon}{3}}P_0\right],\ P_0=\sqrt{\frac{T}{2\pi}} ,
\ee
$0<\epsilon$ is arbitrarily small and ($p$ is the prime)
\be \label{1.2}
\begin{split}
& \zeta_1(s)=\prod_{p\leq P}\sum_{k=0}^\beta\frac{1}{p^{sk}}=\sum_{n<P_0,p\leq P}\frac{1}{n^s}=\sideset{}{'}\sum_{n<P_0}\frac{1}{n^s} , \\
& \zeta_2(s)=\prod_{p\leq P}\sum_{k=0}^\beta\frac{1}{p^{(1-s)k}}=\sum_{n<P_0,p\leq P}\frac{1}{n^{1-s}}=\sideset{}{'}\sum_{n<P_0}\frac{1}{n^{1-s}} , \\
& \zeta_3(s)=\chi(s)\zeta_2(s)
\end{split}
\ee
where (see \cite{3}, p. 16)
\be \label{1.3}
\chi(s)=\pi^{s-\frac 12}\frac{\Gamma\left(\frac{1-s}{2}\right)}{\Gamma\left(\frac s2\right)},\ s\not=2k+1,\ k=0,1,2,\dots \ ,
\ee
and $s=\sigma+it\in\mbb{C}$. We define the function $\tilde{\zeta}(s)$ as follows
\be \label{1.4}
\begin{split}
& \tilde{\zeta}(s)=\tilde{\zeta}(s;P,\beta)=\zeta_1(s)+\zeta_3(s)= \\
& =\sideset{}{'}\sum_{n<P_0}\frac{1}{n^s}+\chi(s)\sideset{}{'}\sum_{n<P_0}\frac{1}{n^{1-s}},\ s\in\mbb{C},\ s\not=2k+1 .
\end{split}
\ee
Since
\bdis
\tilde{\zeta}(1-s)=\sideset{}{'}\sum_{n<P_0}\frac{1}{n^{1-s}}+\chi(1-s)\sideset{}{'}\sum_{n<P_0}\frac{1}{n^s} ,
\edis
and (see \cite{3}, p. 16) $\chi(s)\chi(1-s)=1$, then
\bdis
\tilde{\zeta}(s)=\chi(s)\tilde{\zeta}(1-s),\ s\in\mbb{C},\ s\not= 2k+1 .
\edis

\begin{remark}
The function $\tilde{\zeta}(s)$ obeys the functional equation
\bdis
\tilde{\zeta}(s)=\chi(s)\tilde{\zeta}(1-s)
\edis
hence, the zeros of $\tilde{\zeta}(s)$ either lie on the critical line $\sigma=\frac 12$ or occur in pairs symmetrical about this line.
\end{remark}

\subsection{}

Since (comp. \cite{3}, p. 79)
\bdis
\chi\left(\frac 12+it\right)=e^{-i2\vartheta(t)}
\edis
then from (\ref{1.4}) the formula
\be \label{1.5}
\begin{split}
& e^{i\vartheta(t)}\tilde{\zeta}\left(\frac 12+it\right)=\sideset{}{'}\sum_{n<P_0}\frac{e^{i\{\vth(t)-t\ln n\}}}{\sqrt{n}}+
\sideset{}{'}\sum_{n<P_0}\frac{e^{-i\{\vth(t)-t\ln n\}}}{\sqrt{n}}=\\
& = 2\sideset{}{'}\sum_{n<P_0}\frac{1}{\sqrt{n}}\cos\{\vth(t)-t\ln n\}=Z_1(t;P,\beta)
\end{split}
\ee
follows. We have studied the zeros of $Z_1(t)$, i.e. the zeros of $\tilde{\zeta}(s)$, on the critical line in the paper \cite{2}. Let
\be \label{1.6}
\begin{split}
& D=D(T,H,K)=\{ s:\ \sigma\in [-K,K],\ t\in [T,T+H]\},\\
& K>1,\ T>0,\ H\leq \sqrt{T} .
\end{split}
\ee
In this paper we prove the following theorem.

\begin{mydef1}
\be \label{1.7}
\tilde{\zeta}(s)\not=0,\ s\in D,\ \sigma\not=\frac 12
\ee
for all sufficiently big $T>0$, i.e. for $\tilde{\zeta}(s),\ s\in D,\ T\to\infty$ the analog of the Riemann hypothesis is true.
\end{mydef1}

Let us remind the approximate functional equation of Riemann-Hardy-Littlewood (\cite{3}, p. 69)
\be \label{1.8}
\zeta(s)=\sum_{n\leq t'}\frac{1}{n^s}+\chi(s)\sum_{n\leq t'}\frac{1}{n^{1-s}}+\mcal{O}(t^{-\frac{\sigma}{2}}),\ t'=\sqrt{\frac{t}{2\pi}} ,
\ee
and the Riemann-Siegel formula (comp. (\ref{1.5}))
\be \label{1.9}
\begin{split}
& e^{i\vth(t)}\zf =Z(t)=2\sum_{n\leq t'}\frac{1}{\sqrt{n}}\cos\{\vth(t)-t\ln n\}+\mcal{O}(t^{-\frac 14})= \\
& = 2\sum_{n<P_0}\frac{1}{\sqrt{n}}\cos\{\vth(t)-t\ln n\}+\mcal{O}(T^{-\frac 14})+\mcal{O}(HT^{-\frac 34}),\ t\in [T,T+H] .
\end{split}
\ee
\begin{remark}
The term \emph{the part of the function} $\zeta(s)$ is specified by the comparison of the formulae (\ref{1.4}), (\ref{1.8}). Next, the condition
$H\leq \sqrt{T}$ is related with (\ref{1.9}).
\end{remark}

\section{The formulae for some parts of $\tilde{\zeta}(s)$}

We have (see (\ref{1.2}))
\bdis
\zeta_1(s)=B_1(s)e^{i\psi_1(s)},\ B_1(s)=|\zeta_1(s)|>0,\ \sigma>0 ,
\edis
where
\be \label{2.1}
\begin{split}
& B_1(s)=\prod_{p\leq P}|M_1(p;s,\beta)|,\ \psi_1(s)=\sum_{p\leq P}\arg\{M_1(p;s,\beta)\}, \\
& M_1(p)=\frac{1-Q_1^{\beta+1}}{1-Q_1},\ Q_1=Q_1(p;s)=\frac{1}{p^s},\ |Q_1|=\frac{1}{p^\sigma}<1 ,
\end{split}
\ee
and similarly,
\bdis
\zeta_2(s)=B_2(s)e^{i\psi_2(s)},\ B_2(s)>0,\ \sigma<1 ,
\edis
where
\be \label{2.2}
\begin{split}
& B_2(s)=\prod_{p\leq P}|M_2(p)|,\ \psi_2(s)=\sum_{p\leq P}\arg\{M_2(p)\}, \\
& M_2(p)=\frac{1-Q_2^{\beta+1}}{1-Q_2},\ Q_2=\frac{1}{p^{1-s}},\ |Q_2|=\frac{1}{p^{1-\sigma}}<1 .
\end{split}
\ee
Next, we have (see \cite{3}, pp. 68,79, 329)
\be \label{2.3}
\chi(s)=\left(\frac{t}{2\pi}\right)^{\frac 12-\sigma}e^{-i2\vth(t)}\left\{ 1+\mcal{O}\left(\frac 1t\right)\right\} ,
\ee
i.e.
\bdis
\chi(s)=|\chi(s)|e^{i\psi_3(s)}
\edis
where
\be \label{2.4}
\begin{split}
& |\chi(s)|=\left(\frac{t}{2\pi}\right)^{\frac 12-\sigma}\left\{ 1+\mcal{O}\right\},\ |\chi(s)|>0,\ s\in D, \\
& \psi_3(s)=-2\vth(t)+\mcal{O}\left(\frac 1t\right),\ T\to\infty .
\end{split}
\ee
Consequently, we obtain the following formulae
\be \label{2.5}
\begin{split}
& \tilde{\zeta}(s)=B_1(s)e^{i\psi(s)}+B_2(s)|\chi(s)|e^{i\psi_4(s)} , \\
& \psi_4(s)=\psi_2(s)+\psi_3(s),\ a\in D\cap\{ 0<\sigma<1\},\ T\to\infty .
\end{split}
\ee

\begin{remark}
Let us remind that the formula (\ref{2.3}) is connected with the Stirling's formula for $\ln\Gamma(z),\ z\in\mbb{C}$ to which corresponds arbitrary
fixed strip $-K\leq \sigma\leq K$ (comp. \cite{3}, p. 68).
\end{remark}

\section{The lemmas on $B_1(s),\ B_2(s)$}

\subsection{}

Let
\be \label{3.1}
D_1(\Delta)=\left\{ s:\ \sigma\in \left[\frac 12+\Delta,1-\Delta\right],\ t\in [T,T+H]\right\},\ \Delta\in \left( 0,\frac 14\right) .
\ee
The following lemma holds true.
\begin{lemma}
\be \label{3.2}
\exp\left(-\frac{A}{\Delta}p^{\frac 12-\Delta}\right)<B_1(s)<\exp\left(\frac{A}{\Delta}p^{\frac 12-\Delta}\right),\ s\in D_1(\Delta),\ T\to\infty .
\ee
\end{lemma}

\begin{proof}
We have (see (\ref{2.1}))
\bdis
\begin{split}
& |M_1|=\left| 1-\frac{1}{p^{^{(\beta+1)s}}}\right|\left|1-\frac{1}{p^s}\right|^{-1}=\\
& =\left\{ 1+\frac{1}{p^{2(\beta+1)\sigma}}-\frac{2\cos\{(\beta+1)\vp\}}{p^{(\beta+1)\sigma}}\right\}^{\frac 12}
\left\{ 1+\frac{1}{p^{2\sigma}}-\frac{2\cos\vp}{p^{\sigma}}\right\}^{-\frac 12}=M_{11}M_{12}
\end{split}
\edis
where
\bdis
\vp=t\ln p .
\edis
Next, we have (see (\ref{1.1}))
\be \label{3.3}
\begin{split}
& \ln M_{11}=\frac 12\ln\left\{ 1+\mcal{O}\left(\frac{1}{p^{\frac{\beta}{2}}}\right)\right\}=\mcal{O}\left(\frac{1}{p^{\frac{\beta}{2}}}\right) , \\
& M_{11}=\exp\left\{\mcal{O}\left(\frac{1}{p^{\frac{\beta}{2}}}\right)\right\}
\end{split}
\ee
uniformly for $\Delta\in (0,\frac 14)$, and since
\bdis
\frac 12+\Delta\leq \sigma\leq 1-\Delta,\ \frac{1}{p^{2\sigma}}\leq \frac{1}{p^{1+2\Delta}}<1 ,
\edis
then
\bdis
\begin{split}
& \ln M_{12}=-\frac 12\ln\left( 1-\frac{1}{p^{2\sigma}}\right)-\frac 12\ln\left( 1-\frac{2p^\sigma}{p^{2\sigma}+1}\cos\vp\right)= \\
& =\frac{1}{p^\sigma}\cos\vp+\mcal{O}\left(\frac{1}{p^{2\sigma}}\right) , \\
& M_{12}=\exp\left\{\frac{1}{p^\sigma}\cos\vp+\mcal{O}\left(\frac{1}{p^{2\sigma}}\right)\right\} .
\end{split}
\edis
Hence (see (\ref{2.1})), we have
\be \label{3.4}
B_1(s)=\exp\left\{\sum_{p\leq P}\frac{1}{p^\sigma}\cos\vp+\mcal{O}\left(\sum_{p\leq P}\frac{1}{p^{2\sigma}}\right)\right\},\ s\in D_1(\Delta)
\ee
uniformly for $\Delta\in (0,\frac 14)$. Since
\bdis
\Delta\leq 1-\sigma\leq \frac 12-\Delta ,
\edis
then
\be \label{3.5}
\begin{split}
& \left|\sum_{p\leq P}\frac{1}{p^\sigma}\cos\vp+\mcal{O}\left(\sum_{p\leq P}\frac{1}{p^{2\sigma}}\right)\right|<
A \sum_{p\leq P}\frac{1}{p^\sigma}<\frac{A}{1-\sigma}p^{1-\sigma}< \frac{A}{\Delta}p^{\frac 12-\Delta} ,
\end{split}
\ee
and from this (see (\ref{3.4})) we obtain (\ref{3.2}).
\end{proof}

\subsection{}

The following lemma holds true

\begin{lemma}
\be \label{3.6}
\exp\left(-AP^{1-\Delta}\right)<B_2(s)<\exp\left( AP^{1-\Delta}\right),\ s\in D_1(\Delta),\ T\to\infty
\ee
if the condition
\be \label{3.7}
\Delta\beta>\omega(T)
\ee
is fulfilled, where $\omega(T)$ increases to $\infty$ for $T\to\infty$.

\end{lemma}

\begin{proof}
Since by (\ref{3.7}), (see (\ref{2.2})),
\bdis
(1-\sigma)(\beta+1)\geq \Delta(\beta+1)>\omega(T)
\edis
then putting $1-\sigma=\bar{\sigma}$, we obtain the formula
\be \label{3.8}
B_2(s)=\exp\left\{\sum_{p\leq P}\frac{1}{p^{\bar{\sigma}}}\cos\vp+\mcal{O}\left(\sum_{p\leq P}\frac{1}{p^{2\bar{\sigma}}}\right)\right\} ,
\ee
(similarly to (\ref{3.4})). Since (see (\ref{3.1}), comp. (\ref{3.5}); $\Delta\leq \bar{\sigma}\leq\frac 12-\Delta$)
\bdis
\sum_{p\leq P}\frac{1}{p^{\bar{\sigma}}}+\mcal{O}\left(\sum_{p\leq P}\frac{1}{p^{2\bar{\sigma}}}\right)=
\mcal{O}\left(\sum_{p\leq P}\frac{1}{p^{\bar{\sigma}}}\right)=\mcal{O}\left(\frac{P^{1-\bar{\sigma}}}{1-\bar{\sigma}}\right)=\mcal{O}
(P^{1-\Delta}) ,
\edis
then we obtain (\ref{3.6}) from (\ref{3.8}).
\end{proof}

\begin{remark}
The estimate (\ref{3.6}) is valid in somehow wider domain
\bdis
D_1^+(\Delta)=\left\{ s:\ \sigma\in \left[\frac 12,1-\Delta\right],\ t\in [ T,T+H]\right\} .
\edis
\end{remark}

\section{The function $\tilde{\zeta}(s)$ has no zero in the rectangle $D_1(\Delta_0)$}

First off all (see (\ref{2.4}))

\be \label{4.1}
|\chi(s)|<\frac{A}{P_0^{2\Delta}},\ s\in D_1(\Delta) .
\ee

Next (see (\ref{2.5}), (\ref{3.2}), (\ref{3.6}))

\be \label{4.2}
\begin{split}
& |\tilde{\zeta}(s)|\geq B_1(s)-|\chi(s)|B_2(s)>\exp\left(-\frac{A}{\Delta}P^{\frac 12-\Delta}\right)-\frac{A}{P_0^{2\Delta}}\exp(AP^{1-\Delta})> \\
& > \exp\left(-\frac{A}{\Delta}P^{\frac 12-\Delta}\right)-\frac{A}{P_0^{2\Delta}}\exp(AP)=\\
& = \left\{ 1-\frac{A}{P_0^{2\Delta}}\exp\left( AP+\frac{A}{\Delta}P^{\frac 12-\Delta}\right)\right\}
\exp\left(-\frac{A}{\Delta}P^{\frac 12-\Delta}\right)> \\
& > \left\{ 1-\frac{A}{P_0^{2\Delta}}\exp\left(\frac{2A}{\Delta}P\right)\right\}\exp\left(-\frac{A}{\Delta}P^{\frac 12-\Delta}\right)= \\
& = \left\{ 1-A\exp\left(\frac{2A}{\Delta}P-2\Delta\ln P_0\right)\right\}\exp\left(-\frac{A}{\Delta}P^{\frac 12-\Delta}\right) .
\end{split}
\ee

Since
\bdis
2\Delta\ln P_0-\frac{2A}{\Delta}P=\frac{2\ln P_0}{\Delta}\left( \Delta^2-A\frac{P}{\ln P_0}\right)
\edis
then we put (see (\ref{1.1}))
\be \label{4.3}
\Delta_0=\Delta_0(T,\epsilon)=\left(2A\frac{P}{\ln P_0}\right)^{\frac 12}=\frac{\sqrt{2A}}{(\ln P_0)^{\frac{\epsilon}{2}}} .
\ee

Because (see (\ref{1.1}))
\bdis
\Delta_0\beta > A_1(\ln P_0)^{\frac{\epsilon}{6}} \to\infty,\ T\to\infty
\edis
then the condition (\ref{3.7}) is fulfilled. Hence, we obtain from (\ref{4.2}) by (\ref{1.1}) and (\ref{4.3}) the estimate
\bdis
\begin{split}
& |\tilde{\zeta}(s)|>\frac 12\exp\left(-\frac{A}{\Delta_0}P^{\frac 12-\Delta_0}\right)>\exp\left(-\frac{A}{\Delta_0}P^{\frac 12}\right)=
\exp\left(-\sqrt{\frac A2\ln P_0}\right) , \\
& s\in D_1(\Delta_0),\ T\to\infty .
\end{split}
\edis
Namely, we have the following lemma holds true.

\begin{lemma}
\bdis
|\tilde{\zeta}(s)|> e^{-\sqrt{A\ln P_0}},\ s\in D_1(\Delta_0),\ T\to\infty .
\edis
\end{lemma}

\begin{cor}
\be \label{4.4}
\tilde{\zeta}(s)\not=0,\ s\in D_1(\Delta_0),\ T\to\infty .
\ee
\end{cor}

\section{The function $\tilde{\zeta}(s)$ has no zero in the rectangle $D_2(\Delta_0)$}

Let
\bdis
D_2(\Delta_0)=\{ s:\ \sigma\in [1-\Delta_0,K],\ t\in [T,T+K]\} .
\edis
We remark that the formula (\ref{3.4}) is valid for all $\sigma\in [1-\Delta_0,K]$, see the proof of the Lemma 1. Since in our case (comp. (\ref{3.5}))
\bdis
\left|\sum_{p\leq P}\frac{1}{p^\sigma}\cos\vp+\mcal{O}\left(\sum_{p\leq P}\frac{1}{p^{2\sigma}}\right)\right|<A\sum_{p\leq P}\frac{1}{p^{1-\Delta_0}}
<\frac{A}{\Delta_0}P^{\Delta_0} ,
\edis
then we obtain the estimate (comp. (\ref{3.2}))
\be \label{5.1}
\exp\left(-\frac{A}{\Delta_0}P^{\Delta_0}\right)<B_1(s)<\exp\left(\frac{A}{\Delta_0}P^{\Delta_0}\right)
\ee
for $s\in D_2(\Delta_0),\ T\to\infty$. \\

Next, for $\zeta_2(s)$ we use the formula (see (\ref{1.2}))
\bdis
\zeta_2(s)=\sideset{}{'}\sum_{n<P_0}\frac{1}{n^{1-s}} .
\edis
First of all (see (\ref{1.1}), (\ref{4.3}) and (\ref{1.2}) - the product formula for $\zeta_2(s)$)
\be \label{5.2}
\begin{split}
& \sum_{n<P_0}1=(\beta+1)^{\pi(P)}=\exp\{\pi(P)\ln(\beta+1)\}< \\
& < \exp\left( A(\epsilon)\frac{P}{\ln P}\ln\ln P_0\right)=\exp\left\{ A(\epsilon)(\ln P_0)^{1-\epsilon}\right\}<
\exp(\Delta_0\ln P_0)=P_0^{\Delta_0}
\end{split}
\ee
where we have used the upper estimate of Chebyshev for $\pi(x)$. Next, (see (\ref{1.1}), (\ref{2.4}))
\bdis
|\chi(s)|<\frac{A}{P_0^{2\sigma-1}},\ 1-\Delta_0\leq \sigma\leq K .
\edis
Now:
\begin{itemize}
\item[(A)] in the rectangle
\bdis
D_{21}(\Delta_0)=D_2(\Delta_0)\cup \{ 1-\Delta_0\leq \sigma\leq 1\}
\edis
we have (see (\ref{1.2}), (\ref{5.2}); $1-2\Delta_0\leq 2\sigma-1\leq 1$)
\be \label{5.3}
\begin{split}
& \zeta_3(s)=\chi(s)\zeta_2(s)=\mcal{O}\left(\frac{1}{P_0^{2\sigma-1}}\sideset{}{'}\sum_{n<P_0}\frac{1}{n^{1-\sigma}}\right)=
\mcal{O}\left(\frac{1}{P_0^{2\sigma-1}}\sideset{}{'}\sum_{n<P_0} 1\right)=\\
& = \mcal{O}\left(\frac{1}{P_0^{1-2\Delta_0}}P_0^{\Delta_0}\right)=\mcal{O}\left(\frac{1}{P_0^{1-3\Delta_0}}\right) ,
\end{split}
\ee

\item[(B)] in the rectangle
\bdis
D_{22}(\Delta_0)=D_2(\Delta_0)\cup\{ 1<\sigma\leq K\}
\edis
we have
\be \label{5.4}
\begin{split}
& \zeta_3(s)=\mcal{O}\left(\frac{1}{P_0^{2\sigma-1}}\sum_{n<P_0}\frac{1}{n^{1-\sigma}}\right)=
\mcal{O}\left\{\frac{1}{P_0^\sigma}\sideset{}{'}\sum_{n<P_0}\left(\frac{n}{P_0}\right)^{\sigma-1}\right\} = \\
& = \mcal{O}\left(\frac{1}{P_0^\sigma}\sideset{}{'}\sum_{n<P_0} 1\right)=\mcal{O}\left(\frac{1}{P_0^{1-\Delta_0}}\right) .
\end{split}
\ee
\end{itemize}

Consequently (see (\ref{5.3}), (\ref{5.4})), we have
\be \label{5.5}
\zeta_3(s)=\mcal{O}\left(\frac{1}{P_0^{1-3\Delta_0}}\right),\ s\in D_2(\Delta_0) .
\ee

Since (see (\ref{1.1}), (\ref{4.3}))
\bdis
\frac{A}{\Delta_0}P^{\Delta_0}=\frac{A}{\sqrt{2A_1}}(\ln P_0)^{\frac{\epsilon}{2}}(\ln P_0)^{(1-\epsilon)\Delta_0}<(\ln P_0)^{\frac{2\epsilon}{3}}, \
T\to\infty ,
\edis
then (see (\ref{2.5}), (\ref{5.1}), (\ref{5.5})) we obtain in the domain $D_2(\Delta_0)$
\bdis
\begin{split}
& |\tilde{\zeta}(s)|\geq B_1(s)-|\zeta_3(s)|>\exp\left(-\frac{A}{\Delta_0}P^{\Delta_0}\right)-\frac{A}{P_0^{1-3\Delta_0}}= \\
& = \left\{ 1-\exp\left[\frac{A}{\Delta_0}P^{\Delta_0}-(1-3\Delta_0)\ln P_0+\ln A\right]\right\}\exp\left(-\frac{A}{\Delta_0}P^{\Delta_0}\right)> \\
& > \frac 12\exp\left[-(\ln P_0)^{\frac{2\epsilon}{3}}\right]>\exp\left[-(\ln P_0)^{\epsilon}\right],\ T\to\infty ,
\end{split}
\edis
i.e. the following lemma holds true.

\begin{lemma}
\bdis
|\tilde{\zeta}(s)|>e^{-(\ln P_0)^\epsilon},\ s\in D_2(\Delta_0),\ T\to\infty .
\edis
\end{lemma}

\begin{cor}
\be \label{5.6}
\tilde{\zeta}(s)\not=0,\ s\in D_2(\Delta_0),\ T\to\infty .
\ee
\end{cor}

\section{Lemma on the difference of logarithms}

Let

\be \label{6.1}
D_3(\Delta_0)=\left\{ s:\ \sigma\in \left(\left. \frac 12,\frac 12+\Delta_0\right.\right],\ t\in [T,T+H]\right\}
\ee
where
\bdis
\sigma=\frac 12+\delta,\ \delta\in (0,\Delta_0) .
\edis
The following lemma holds true.

\begin{lemma}
\be \label{6.2}
\ln B_1(s)-\ln B_2(s)=\mcal{O}\left\{\delta (\ln P_0)^{\frac{1-\epsilon}{2}}\right\},\ s\in D_3(\Delta_0),\ T\to\infty .
\ee
\end{lemma}

\begin{proof}
We have (see (\ref{2.1}), (\ref{2.2}))
\be \label{6.3}
\ln B_1(s)-\ln B_2(s)=Y_1+Y_2
\ee
where ($|z|=|\bar{z}|$)
\bdis
\begin{split}
& Y_1=\sum_{p\leq P}\left\{\ln\left| 1-\frac{p^{-i(\beta+1)t}}{p^{(\beta+1)(\frac 12+\delta)}}\right|-
\ln\left| 1-\frac{p^{-i(\beta+1)t}}{p^{(\beta+1)(\frac 12-\delta)}}\right|\right\} , \\
& Y_2=\sum_{p\leq P}\left\{ \ln\left| 1-\frac{p^{it}}{p^{\frac 12-\delta}}\right|-\ln\left| 1-\frac{p^{it}}{p^{\frac 12+\delta}}\right|\right\} .
\end{split}
\edis
Let
\bdis
x=\frac{1}{p^\sigma},\ x\in \left[\frac{p^{-\delta}}{\sqrt{p}},\frac{p^{\delta}}{\sqrt{p}}\right] .
\edis
It is clear that (see (\ref{1.1}), (\ref{4.3}))
\bdis
\begin{split}
& \delta\ln p=\mcal{O}(\Delta_0\ln P)=\mcal{O}\left\{\frac{\ln\ln P_0}{(\ln P_0)^{\frac{\epsilon}{2}}}\right\}\to 0,\ T\to\infty , \\
& \frac{p^{\delta}-p^{-\delta}}{\sqrt{p}}=\mcal{O}\left(\delta\frac{\ln P_0}{\sqrt{p}}\right) .
\end{split}
\edis
Next, by the mean-value theorem
\bdis
\begin{split}
& \ln\left| 1-\frac{p^{-i(\beta+1)t}}{p^{(\beta+1)(\frac 12+\delta)}}\right|-\ln\left| 1-\frac{p^{-i(\beta+1)t}}{p^{(\beta+1)(\frac 12-\delta)}}\right| = \\
& =\left.\frac{p^{-\delta}-p^\delta}{\sqrt{p}}\frac{{\rm d}}{{\rm d}x}\left\{\ln\left| 1-x^{\beta+1}p^{-i(\beta+1)t}\right|\right\}\right|_{x=x_1} , \\
& x_1=\frac{1}{p^c},\ c\in \left(\frac 12-\delta,\frac 12+\delta\right) .
\end{split}
\edis
Since ($\vp=t\ln p$)
\bdis
\ln\left| 1-x^{\beta+1}p^{-i(\beta+1)t}\right|=\frac 12\ln\left( 1+x^{2\beta+2}-2x^{\beta+1}\cos\{(\beta+1)\vp\}\right) ,
\edis
then (see (\ref{1.1}))
\bdis
\begin{split}
& \frac{{\rm d}}{{\rm d}x}\ln\left| 1-x^{\beta+1}p^{-i(\beta+1)t}\right|= \\
& = \frac 12\frac{(2\beta+2)x^{2\beta+1}-2(\beta+1)x^\beta\cos\{(\beta+1)\vp\}}{1+x^{2\beta+2}-2x^{\beta+1}\cos\{(\beta+1)\vp\}}=
\mcal{O}\left(\frac{\beta}{p^{\sigma\beta}}\right)=\mcal{O}\left(\frac{\beta}{p^{\frac \beta 2}}\right)
\end{split}
\edis
where $\sigma>\frac 12$, and consequently
\be \label{6.4}
Y_1=\mcal{O}\left(\delta\sum_{p\leq P}\frac{\beta}{p^{\frac \beta 3}}\frac{\ln p}{\sqrt{p}}\right)=
\mcal{O}\left(\delta\frac{\beta}{2^{\frac \beta 6}}\sum_{p\leq P}\frac{1}{p^{\frac \beta 6}}\right)=\mcal{O}(\delta) .
\ee
Similarly, we obtain in the case $Y_2$
\bdis
\ln\left| 1-\frac{p^{it}}{p^{\frac 12-\delta}}\right|-\ln\left| 1-\frac{p^{it}}{p^{\frac 12+\delta}}\right|=
\left.\frac{p^{-\delta}-p^{\delta}}{2\sqrt{p}}\frac{{\rm d}}{{\rm d}x}\ln\left( 1+x^2-2x\cos\vp\right)\right|_{x=x_2}
\edis
where
\bdis
\frac{{\rm d}}{{\rm d}x}\ln\left( 1+x^2-2x\cos\vp\right)=\frac{2x-2\cos\vp}{1+x^2-2x\cos\vp}=\mcal{O}(1) ,
\edis
because
\bdis
1+x^2-2x\cos\vp\geq (1-x)^2>\left( 1-2^{-\frac 12+\delta}\right)^2>\left( 1-2^{-\frac 13}\right)^2 > 0 .
\edis
Thus, we have (see (\ref{1.1})
\be \label{6.5}
Y_2=\mcal{O}\left(\delta\sum_{p\leq P}\frac{\ln p}{\sqrt{p}}\right)=\mcal{O}(\delta\sqrt{P})=\mcal{O}\left\{\delta(\ln P_0)^{\frac{1-\epsilon}{2}}\right\} ,
\ee
(the Abel's transformation was used, comp. \cite{2}, (31), (33)). Now, (\ref{6.2}) follows from (\ref{6.3}) by (\ref{6.4}), (\ref{6.5}).
\end{proof}

\section{More accurate formula for $\ln|\chi(s)|$}

The following lemma holds true.

\begin{lemma}
\be \label{7.1}
\ln|\chi(s)|=-\left(\sigma-\frac 12\right)\ln\frac{t}{2\pi}+\mcal{O}\left(\frac{2\sigma-1}{t}\right),\ s\in D_3(\Delta_0),\ T\to\infty .
\ee
\end{lemma}

\begin{proof}
Since (see (\ref{1.3}))
\bdis
|\chi(s)|=\pi^{\sigma-\frac 12}
\left|\frac{\Gamma\left(\frac{1-\sigma}{2}-i\frac t2\right)}{\Gamma\left(\frac{\sigma}{2}+i\frac t2\right)}\right|=\pi^{\sigma-\frac 12}G_1(\sigma,t) ,
\edis
where
\bdis
G_1(\sigma,t)>0,\ s\in D_3(\Delta_0),\ T\to\infty,
\edis
then
\be \label{7.2}
\ln|\chi(s)|=\left(\sigma-\frac 12\right)\ln \pi+\ln G_1(\sigma,t)=\left(\sigma-\frac 12\right)\ln\pi+G_2(\sigma,t),
\ee
and $G_2(\sigma,t)$ is the analytic function of the real variable $\sigma$ for arbitrary fixed $t$ if $s\in D_3(\Delta_0),\ T\to\infty$. Since
\bdis
\left|\chi\left(\frac 12+it\right)\right|=1
\edis
then
\bdis
G_2\left(\frac 12,t\right)=0,
\edis
and
\bdis
G_2(\sigma,t)=\left(\sigma-\frac 12\right)G_3(\sigma,t) .
\edis
Now, (see (\ref{7.2}))
\be \label{7.3}
\ln|\chi(s)|=\left(\sigma-\frac 12\right)\{\ln\pi+G_3(\sigma,t)\},\ s\in D_3(\Delta_0),\ T\to\infty .
\ee
In the case (\ref{2.4}) we have
\be \label{7.4}
\ln|\chi(s)|=-\left(\sigma-\frac 12\right)\ln\frac{t}{2\pi}+G_4(\sigma,t),\ G_4(\sigma,t)=\mcal{O}\left(\frac 1t\right) ,
\ee
under the conditions (\ref{7.3}) where $G_4(\sigma,t)$ is the analytic function of the real variable $\sigma$. Since
\bdis
G_4\left(\frac 12,t\right)=0
\edis
(see (\ref{7.3}), (\ref{7.4})) then
\be \label{7.5}
G_4(\sigma,t)=\left(\sigma-\frac 12\right)G_5(\sigma,t) .
\ee
Next, by (\ref{7.5}), the orders of the functions
\bdis
G_4(\sigma,t),\ G_5(\sigma,t),\ s\in D_3(\Delta_0),\ T\to\infty
\edis
in the variable $t$ are equal, i.e. (see (\ref{7.4}))
\be \label{7.6}
G_5(\sigma,t)=\mcal{O}\left(\frac 1t\right) .
\ee
Now, the formula (\ref{7.1}) follows from (\ref{7.4}) by (\ref{7.5}), (\ref{7.6}).
\end{proof}

\begin{remark}
The formula (\ref{7.1}) can be proved directly, of course.
\end{remark}

\section{Proof of the Theorem}

\subsection{}

Let

\bdis
\ln\Lambda(s)=\ln B_1(s)-\ln B_2(s)-\ln|\chi(s)|,\ s\in D_3(\Delta_0),\ T\to\infty .
\edis
Since (see (\ref{1.1}), (\ref{7.1}), $\sigma=\frac 12+\delta$)
\be \label{8.1}
\ln|\chi(s)|=-\delta\ln\frac{t}{2\pi}+\mcal{O}\left(\frac{\delta}{t}\right)=-2\delta\ln P_0+\mcal{O}\left(\frac{\delta H}{T}\right)+
\mcal{O}\left(\frac{\delta}{T}\right)
\ee
then we obtain (see (\ref{1.6}), (\ref{6.2}), (\ref{8.1}))
\be \label{8.2}
\begin{split}
& \ln\Lambda(s)=2\delta\ln P_0+\mcal{O}\left\{\delta(\ln P_0)^{\frac{1-\epsilon}{2}}\right\}+\mcal{O}\left(\frac{\delta}{\sqrt{T}}\right)=\\
& = \delta\left[ 2\ln P_0+\mcal{O}\left\{ (\ln P_0)^{\frac{1-\epsilon}{2}}\right\}+\mcal{O}\left(\frac{1}{\sqrt{T}}\right)\right]>
\delta\ln P_0>0 .
\end{split}
\ee
Consequently,
\bdis
\Lambda(s)>1,\ s\in D_3(\Delta_0),\ T\to\infty ,
\edis
and (see (\ref{2.5}))
\be \label{8.3}
\begin{split}
& |\tilde{\zeta}(s)|\geq B_1(s)-|\chi(s)|B_2(s)=|\chi(s)|B_2(s)\left(\frac{B_1(s)}{|\chi(s)|B_2(s)}-1\right)=\\
& = |\chi(s)|B_2(s)[\Lambda(s)-1]>0,\ s\in D_3(\Delta_0),\ T\to\infty .
\end{split}
\ee
(The inequality $B_2(s)>0,\ s\in D_3(\Delta_0),\ T\to\infty$ follows from Remark 4.) Now, by (\ref{4.4}), (\ref{5.6}), (\ref{8.3}) and from Remark 1, we have
(\ref{1.7}).

\subsection{}

As an addition to the Theorem we obtain an lower estimate for $|\tilde{\zeta}(s)|,\ s\in D_3(\Delta_0)$. Namely, we have

\begin{lemma}
\be \label{8.4}
|\tilde{\zeta}(s)|>\frac{1}{P_0}\sinh\left(\frac{\delta}{2}\ln P_0\right),\ s\in D_3(\Delta_0),\ T\to\infty,\ \delta\in (0,\Delta_0) .
\ee
\end{lemma}

\begin{proof}
Since (see (\ref{1.1}), (\ref{3.6}), Remark 4 and (\ref{6.1}))
\bdis
B_2(s)>\exp\left(-AP^{1-\Delta_0}\right)>\exp(-AP)=\exp\left\{ -A(\ln P_0)^{1-\epsilon}\right\} ,
\edis
and (see (\ref{8.1}), (\ref{8.2}))
\bdis
|\chi(s)|>P_0^{-(2+\epsilon)\delta},\ \Lambda(s)>P_0^\delta
\edis
then (see (\ref{8.3}))
\bdis
\begin{split}
& |\tilde{\zeta}(s)|>\frac{\exp\left\{-A(\ln P_0)^{1-\epsilon}\right\}}{P_0^{(2+\epsilon)\delta}}\left( P_0^{\delta}-1\right)> \\
& > 2\frac{\exp\left\{-A(\ln P_0)^{1-\epsilon}\right\}}{P_0^{(\frac 32+\epsilon)\Delta_0}}\sinh\left(\frac{\delta}{2}\ln P_0\right)>
\frac{1}{P_0}\sinh\left(\frac{\delta}{2}\ln P_0\right) ,
\end{split}
\edis
i.e. (\ref{8.4}).
\end{proof}

\section{A new property of the functions $\pi(x),\ R(x)$}

Let
\bdis
\begin{split}
& D(\Delta_0)=D(\Delta_0,T,H,K)= \\
& = \left\{ s:\ \sigma\in \left[\frac 12+\Delta_0,K\right],\ t\in [T,T+H]\right\},\ \Delta_0=\frac{A}{(\ln P_0)^{\frac{\epsilon}{2}}} .
\end{split}
\edis
We can prove the following

\begin{mydef61}
\be \label{9.1}
\ln\tilde{\zeta}(s)=s\int_2^P\frac{\pi(x)}{x(x^s-1)}{\rm d}x-\pi(P)\ln\left( 1-\frac{1}{P^s}\right)+\mcal{O}\left(e^{-A\beta}\right) ,
\ee
where $s\in D(\Delta_0),\ T\to\infty$ and $\mcal{O}\left(e^{-A\beta}\right)$ is the estimate of
\bdis
\sum_{p\leq P}\ln\left( 1-\frac{1}{p^{s(\beta+1)}}\right)+\ln\left( 1+\frac{\chi(s)\zeta(s)}{\zeta_1(s)}\right)=\Omega_1(s;P,\beta) .
\edis
\end{mydef61}

Since
\bdis
\pi(x)=\int_0^x\frac{{\rm d}x}{\ln v}+R(x)=U(x)+R(x) ,
\edis
then we obtain from (\ref{9.1})

\begin{mydef62}
\bdis
\begin{split}
& \ln\tilde{\zeta}(s)=s\int_2^P\frac{R(x)}{x(x^s-1)}{\rm d}x-R(P)\ln\left(1-\frac{1}{P^s}\right)+\mcal{O}\left(e^{-A\beta}\right),\\
& s\in D(\Delta_0),\ T\to\infty
\end{split}
\edis
where $\mcal{O}\left(e^{-A\beta}\right)$ is the estimate of
\bdis
\Omega_1(s;P,\beta)+\Omega_2(s;P),
\edis
and
\bdis
\begin{split}
& \Omega_2(s)=-s\int_2^P\frac{{\rm d}x}{\ln x}\int_2^x\frac{{\rm d}v}{v(v^s-1)}-U(P)\ln\left( 1-\frac{1}{2^s}\right)= \\
& =\sum_{n=0}^\infty \frac{1}{n+1}\int_2^P\frac{{\rm d}x}{x^{(n+1)s}\ln x}=\mcal{O}\left(\frac{\sqrt{\ln T}}{T}\right) .
\end{split}
\edis
\end{mydef62}

Thus, the following properties of the functions $\pi(x),\ R(x)$ holds true:
\begin{itemize}
\item[(A)] the function $\pi(x),\ x\in[2,P]$ is the solution of the integral equation (for every fixed $s\in D(\Delta_0)$)
\be \label{9.2}
\ln\tilde{\zeta}(s)=s\int_2^P\frac{\Phi(x)}{x(x^s-1)}{\rm d}x-\Phi(P)\ln\left( 1-\frac{1}{P^s}\right)+\Omega_1(s) ,
\ee
\item[(B)] the function $R(x),\ x\in [2,P]$ is the solution of the perturbed integral equation
\be \label{9.3}
\ln\tilde{\zeta}(s)=s\int_2^P\frac{\Phi(x)}{x(x^s-1)}{\rm d}x-\Phi(P)\ln\left( 1-\frac{1}{p^s}\right)+\Omega_1(s)+\Omega_2(s) .
\ee
\end{itemize}

\begin{remark}
Hence, we have a new property of the functions $\pi(x)$ and $R(x)$: these functions are to solutions of the integral equations (\ref{9.2}) and
(\ref{9.3}), respectively and the mentioned integral equations are close each to other. This property of $\pi(x)$ and $R(x)$ is fully missing in the
theory of $\pi(x),\ R(x)$ based on the Riemann zeta-function.
\end{remark}

Let us remind that the behaviour of the functions $\pi(x),\ R(x)$ is essentially different, $\pi(x)\sim \frac{x}{\ln x},\ x\to\infty$, and
$R(x),\ x\to\infty$ infinitely many times alternates its sign (Littlewood, 1914). \\

\thanks{I would like to thank Michal Demetrian for helping me with the electronic version of this work.}

\end{document}